# Situating "Ethics in Mathematics" as a Philosophy of Mathematics Ethics Education

## Dennis Müller
RWTH Aachen University, Germany
dennis.mueller3@rwth-aachen.de

## Abstract

In this paper, we situate the educational movement of "Ethics in Mathematics," as outlined by the Cambridge University Ethics in Mathematics Project, in the wider area of mathematics ethics education. By focusing on the core message coming out of Ethics in Mathematics, its target group, and educational philosophy, we set it into relation with "Mathematics for Social Justice" and Paul Ernest's recent work on ethics of mathematics. We conclude that, although both Ethics in Mathematics and Mathematics for Social Justice appear antagonistic at first glance, they can be understood as complementary rather than competing educational strategies.

**Keywords:** Ethics in Mathematics, Mathematics for Social Justice, practical ethics,

# Introduction

Ethics and mathematics are often perceived as two entirely different disciplines: one dealing with the questions of leading a moral life and the other dealing with abstract truths (Ernest 2020c). Particularly, pure mathematics is often seen as ethically pure and neutral, posing a challenge to those engaging in educational work on ethics of mathematics (Ernest 2020a). Few philosophers study the making of mathematical knowledge with a focus on ethics (Rittberg 2020), and ethics courses specifically designed for mathematicians can be found in even fewer universities. Understandably so: Questions about the ethics of mathematics, like all ethical questions, entail explosive potential. They can be challenging and uncomfortable. Members of the community had a first-hand experience when Colin Rittberg spoke about the notion of epistemic injustice at the 2nd Workshop on Ethics in Mathematics in Cambridge (Rittberg et al. 4/3/2019). Building on anecdotal evidence, Rittberg presented the case study of Olivia Caramello who tried to publish a paper about a theorem whose proof had never been formally published. Her initial publication was initially denied since the result was supposed to be "common knowledge" among a small group of experts. She later on successfully published her work (Rittberg et al. 2018). Some workshop participants familiar with the Caramello case left the room before Rittberg began his talk. Rittberg described these events as a schism in the conference:

> "[S]ome mathematicians got quite upset when I gave this case study. […] It was a conference that had a schism and was by far the most exciting conference I ever went to." (Rittberg 10/25/2021)

Rittberg's notion of epistemic injustice does not go unchallenged in the mathematics education literature, and some argue in favor of calling it epistemic exclusion (Baldino and Cabral 2021), a term which possibly would have led to less controversy among the conference participants. At the conference, Rittberg's idea of epistemic injustice clashed with a commonly found understanding of mathematics as an object rather than social practice. Pitsili-Chatzli (2021) outlined how through different methods of discourse, these views can be antagonistic when they try to fix a specific meaning of mathematics and its elements of proof and formality. Rittberg's presentation critically examined the role of culture, power, and ideology in pure mathematics and challenged its practitioners directly.

However, the schism did not only occur through the explosive nature of the topic but also because it set the stage for a heated debate: What is the right way to understand the ethical issues of mathematics, its teaching, and practice? Some leaned towards what we describe as "Ethics in Mathematics," while others proclaimed "Mathematics for Social Justice" as the necessary point of view. In this paper, we try to situate Ethics in Mathematics in this debate. We first outline three different educational movements: Ethics in Mathematics (EiM), as developed by the Cambridge University Ethics in Mathematics Project, Mathematics for Social Justice, and Paul Ernest's recent work on ethics of mathematics. After that, we further contrast the educational philosophy behind Ethics in Mathematics with that behind Mathematics for Social Justice. We conclude that Ethics in Mathematics and Mathematics for Social Justice can be seen as complementary rather than competing educational strategies.

# Overview

## Ethics in Mathematics

The contrast between, at first glance, ethically neutral and aesthetically pure mathematical objects and mathematics as human activity has again come into focus and is being challenged in the small but growing research area called "Ethics in Mathematics." Mathematical practice is characterized by remarkable freedom: mathematicians create the definitions with which they later continue to work. Even if proofs, understood line by line, are often a stringing together of trivial statements and their consequences, the overall construct of the proof, from its approach to the mathematical tools used and its logical structuring, is nevertheless characterized by fundamental freedom (Nickel 2005). Above all, mathematical practice is thus a particular form of human action. Ethics in Mathematics takes on this idea and presupposes that human actions can never be ethically neutral, and, hence, there are ethical consequences in every aspect of mathematics, from its definitions and theorems to its teaching and its applications, that should be studied and taught. In the following, we outline "Ethics in Mathematics" as described by the Cambridge University Ethics in Mathematics Society (CUEIMS 2016) and the Cambridge University Ethics in Mathematics Project (CUEiMP 2018).[1]

## Four Levels of Ethical Engagement

Much of the early works on Ethics in Mathematics deal with the question of "Is there ethics in mathematics?" (Müller 2018). The mathematician's possible negative reply to this question is taken as a challenge from which Chiodo and Bursill-Hall (2018) build a level-based system for a mathematician's ethical engagement:

- Level 0: Believing there is no ethics in mathematics
- Level 1: Realising there are ethical issues inherent in mathematics

---

[1] Other courses and projects have taken on a slightly different position on Ethics in Mathematics.

- Level 2: Doing something: speaking out to other mathematicians
- Level 3: Taking a seat at the tables of power
- Level 4: Calling out the bad mathematics of others

The ordering of these levels is deeply related to common definitions of good and bad mathematics. Mathematicians commonly define "good" mathematics using the logical "OR" operator: mathematics is good if it satisfies at least one of N aspects. Tao's reflections on good mathematics are representative of this view. He lists 21 different aspects of good mathematics ranging from good problem-solving and technique, exposition, taste, and pedagogy, to deep, strong, and intuitive mathematics (Tao 2007). Mathematics can thus be judged on utilitarian grounds when it helps to build other mathematics and on intrinsic values, including aesthetics. Most importantly, however, good mathematics is mostly judged within mathematics, without reference to external societal and ethical factors.

But while the common notions of good mathematics are intra-mathematical, those of bad mathematics often not only reference logical mistakes but specifically include negative societal and ethical consequences. For example, Mouhot outlined a typology of bad mathematics that included following a wrong argument, building on false axioms, or practicing math-science washing to boost an opinion by making it look more scientific and mathematical (Mouhot 4/4/2019). Through the latter, practitioners can give mathematics "excessive authority and direct it in ways that [can] cause harm and exploit others" (Chiodo and Müller 2018, p. 6). This leads to a first problem: what mathematicians consider good might only be good from a limited perspective, while it can appear bad by changing one's point of view and by considering extra-mathematical consequences and motivational aspects. Consider, for example, an explicit construction or even hypothetical proof of fast integer factorization. It would undeniably provide deep insight into mathematics, and hence, it would correctly be identified as good mathematics according to Tao's definition. However, the societal consequences of it could be deleterious because much of our modern digital security builds on this very problem being hard to solve.

This common definitional approach to good and bad mathematics leads to a second related problem for anyone considering the ethics of mathematics. It raises a communicational barrier between the ethicist and the mathematician. Our common-sense understanding of good and bad is mutually exclusive: something is either good or bad, but it is not both. However, if you define "good" using the logical "OR" operator and exclude extra-mathematical consequences and motivational aspects, this implies that the definition of "bad" must be "not one of the 21 aspects is satisfied." But this will likely never be satisfied for correct higher mathematics. From an ethics perspective, a much better definition of good mathematics would therefore be split up into "ethically good" and "mathematically good," i.e., "mathematics is good if it satisfies at least one of N aspects, and if at the same time, it does not do one of these M ethically bad aspects." The standard definitions of bad mathematics do not adequately address this issue. Standing on their own, Mouhot's and Tao's definitions appear sound at first, but when judged in combination, they are insufficient to talk about the ethics of mathematics. This proves to be an immense communicative challenge that implicitly finds much attention in Ethics in Mathematics. Its difficulty is a critical aspect to why calling out bad mathematics is identified as the highest level of ethical awareness for mathematicians.

A primary desideratum of any work in ethics of mathematics is currently the lack of data. There are few qualitative studies on research mathematics, mathematicians, and their understanding of ethics. Existing studies of research mathematicians are limited to questions about basic philosophical assumptions (Sriraman 2004), the art of doing mathematics (Heintz 2000b, 2000a; Kiesow 2016), including the impact of computational and AI-driven maths (Sørensen 2016; Ebert 2016) and about the understanding of numbers (Heintz 2007). Existing studies about ethical mathematics understanding focus on high school students and non-mathematics majors (Stephan et al. 2021; Demattè 2021). To the author's knowledge, questions about the ethical self-understanding of the discipline have neither been explored through quantitative nor qualitative interview-based studies. Only a little work has gone into why beliefs about the certainty of mathematical knowledge and ethical neutrality seem to be so widespread among students and practitioners alike (Ernest 2016a). Recent research has analyzed the public attitudes towards artificial intelligence, but studies about expert practitioners' attitudes towards AI and its ethics are still missing (Ikkatai et al. 2022). Further research is urgently needed to shed light on the apparently common practitioner's belief that "there is no ethics in mathematics," which Chiodo and Bursill-Hall (2018) described in the first Ethics in Mathematics Discussion paper.

## Core Message: Preventing Harm

Müller et al. (2021) outlined that until a few years ago, the normative studies that mathematicians undertook about mathematics consisted mainly of many investigations of individual sub-fields. Over the years, individual mathematical disciplines came into the spotlight of ethical attention, including statistics, mathematical physics, financial mathematics, cryptography, and today's widely popular machine learning. Most of these ethical analyses have one thing in common. Rather often, there was quite a loud bang beforehand - sometimes proverbially, as with the atomic bombs dropped on Japan at the end of the Second World War, sometimes metaphorically, as in the global financial crisis of 2008. But in contrast to these events, many of the involved mathematical disciplines were very abstract in their origins and were widely considered ethically neutral until their applications became apparent. Indeed, the line between applied and pure mathematics can be "irritatingly vague" (Hersh 1990, p. 22) and very much depends on the mathematical questions asked (Brown 2008). Thus, even during World War II and only a few years before the development of atomic bombs, Hardy still referred to quantum mechanics as an abstract and useless science (Hardy 1940). Number theory, too, long had the whiff of the ethically pure about it. Beginning with Plato and ending with the most abstract considerations in algebraic number theory, for a long time, these considerations had little influence on the everyday life of most people. However, at the latest, just how big the ethical questions are that number theory and cryptography have imposed on our modern digital society have become apparent after the many recent global surveillance disclosures.

Hence, it is not surprising that both the Cambridge University Ethics in Mathematics Project and the Cambridge University Ethics in Mathematics Society pursue "preventing harm" as their declared objective:

> "It is […] obvious that high-powered and sophisticated mathematics is ubiquitous in modern technology, finance, the nation's infrastructure and defense, and social media. And pretty much anywhere else in the 21st century. Much of mathematics can be used for good – but simply put: it can also be a tool for harm. We think it is clear that mathematical research and practice may take its practitioners to deep and professionally specific ethical issues." (CUEiMP 2018)

The old view of ethics of mathematics as sub-field investigations is consistent with a problematic approach found in modern risk analysis. It frequently lacks focus on structural risks and puts too much attention on the known unknowns: we ask ourselves if this one technology, one drug, or one fertilizer is safe, and we too rarely ask how individual technologies and methods interact structurally (Jasanoff 2016). Ethics in Mathematics attempts to avoid this trap by abstracting from subdisciplines and studying the ethical aspects valid for all mathematical practitioners (CUEiMP 2018). By its founders, this new approach is perceived as sorely needed (Chiodo and Clifton 2019).

### Practical Ethics for the Working Mathematician

Ethics in Mathematics specifically targets working mathematicians and tries to prepare them to engage in societal discourses about their field. For instance, the close connection between civilian and military research in mathematics (e.g., Gruber 2017, 2018) was put into the spotlight during the 2017 to 2019 CUEiMS Ethics in Mathematics seminars. Through a series of guest speakers ranging from whistle-blowers, such as Edward Snowden, William Binney, and Julian Assange, to mathematicians and computer scientists, including Richard Stallman, Tanja Lange, and Martin Hellman (CUEIMS 2016), the society complemented the regular lectures to hear from frontline practitioners about mathematical responsibility.

Since its early years, both the Cambridge Project and Society have also intently addressed undergraduate mathematicians without any formal training in philosophy while maintaining a certain degree of neutrality in their educational approach. Unlike Mathematics for Social Justice, neither the Project nor the Society has actively taken on a specific theory from education. Their understanding of neutrality went so far that it was suggested to address ethics issues in mathematics within the framework of a student's definition of pure mathematics instead of challenging a possibly misleading or false definition itself (Müller 2018). Both project and society are interdisciplinary, if not transdisciplinary, in their educational approach, inviting speakers from various backgrounds to their seminar series and conferences (CUEiMP 2018; CUEIMS 2016). Ernest writes about the Cambridge University Ethics in Mathematics Project that

> "[o]verall, this project can be seen as a valuable and exemplary first attempt to plug the hole left by the widespread neglect of ethics for mathematicians. It is primarily practice-driven and does not include any elements of theoretical ethics. This means that good and bad, vice and virtue are not analyzed. It is assumed that bad or harmful applications can be identified as such through everyday ethical common-sense. The organizers are probably right in making this judgment call with regard to developing an undergraduate course." (Ernest 2021, p. 20)

The description for its eight lectures can be found in the term cards of the student society (CUEIMS 2016). They are titled as follows:

1. An introduction to ethics in mathematics and why it is important
2. Financial mathematics and modeling
3. Cryptography, surveillance, and privacy
4. Fairness and impartiality in algorithms and AI
5. Regulation, accountability, and the law
6. Understanding the behavior of the mathematical community
7. Psychology 101: How to survive as a mathematician at work

8. Looking into the future, what more can mathematicians do?

Lecture 7, as described in the Michaelmas 2021 term card, outlines that as future mathematicians, students will work in complex work environments, subject to the "usual workplace interactions, issues, conflicts, and dangers that may arise in other professions." It tries to prepare the students by making them understand their strengths and weaknesses, building on the cultural analysis of the mathematical community in lecture 6. Both lectures can also be understood in what Hersh calls "front mathematics" (Hersh 1991). The formal, precise, ordered, and abstract positions that the field projects to the outside are not necessarily going to protect the individual mathematician in the workplace, as it merely protects the mathematical endeavor as a whole. As such, the lecture series is concerned with making mathematicians perceive the ethical relevance of their mathematics and tries to challenge some conservative tendencies in the field. As Lovász put it, "mathematicians are conservative people; I don't mean that we are right-wing, but we don't push for changes: We are reluctant to spend time on anything other than trying to prove that P != NP." (Lovász 1998, p. 33).

## Relationship to Other Courses

The history of mathematics and history of science classes at the faculty of mathematics at Cambridge University are running in parallel to the ethics seminars. They address the time-varying normative power ascribed to and coming from mathematics through its uses in the natural sciences, including but not limited to changing notions of the relationship between mathematical models and the physical world, and different modes of mathematical reasoning in the biological sciences. In its applications, mathematics has not only a descriptive but also a normative function in many areas of our lives, as it is often used to create and introduce new economic and social rules (Nickel 2013). In these areas, the ethical discourse usually does not end with their mathematization but often raises new ethical questions through its very introduction (Nickel 2014). The introduction of a new mode of discourse extends to the use of mathematics within the natural sciences and engineering. In the former, it is the relationship between mathematical models and nature that introduced new normative elements into society. In the latter, it happens through the consequences coming from mathematical technologies.

The history lectures particularly stress the complexity of the historian's work and the difficulty it takes to understand the slowly changing notions of descriptive and normative power related to mathematics resting on Owen Gingerich's work, among others. Over a 35 year spanning project, Gingerich showed that Copernicus's De revolutionibus orbium coelestium sections on calculations were widely annotated and read by mathematicians and astronomers (Gingerich 1992, 2002, 1993, 2004). Gingerich found severely annotated sections on calculations, but the section on cosmology was lacking marginalia, which suggested to him that it was mostly ignored by contemporary mathematicians and astronomers. Copernicus was likely read as a calculational improvement over Ptolemy's Almagest. The shift that mathematical models should be representative of the physical world had neither reached religious, societal, nor philosophical consensus immediately after Copernicus's publication. The debates went on to Galileo, Kepler, and Newton, who in increasing detail introduced a new understanding of the relationship between mathematical models and the physical world.

Galileo famously exclaimed that the book of nature was written in the language of mathematics and understood through inductive reasoning using experimental knowledge, breaking with the Aristotelean tradition of deductive reasoning from book knowledge. But while Galileo considered repeatable physical experiments from which he inferred the physical laws written in mathematics, Darwin could not build on experiments but merely on his observations of nature. He had to infer the meaning based on what he could observe in nature and combined it with probabilistic reasoning to arrive at his theory of evolution. This subtle difference of Darwin's use of probabilistic thinking in his Origins of Species which is now found in many applications of probability and statistics, was not covered in the ethics seminar, but those interested in it could hear about it in the accompanying history classes. Thus, two of the most important examples of how we see mathematics today were not directly part of the Ethics in Mathematics seminar.

## Mathematics for Social Justice

In contrast, a different approach can be found in the many works on Mathematics for Social Justice, which (for most authors) finds its roots in the more extensive critical pedagogy literature. Critical pedagogy builds on the intimate relationships between culture, power, public discourse, and democracy and their connections to education. Questions of race, class, and gender play a central role in this approach. Therefore, its proponents argue that critical scholars and educators need to take a position that promotes human freedom and social justice (e.g., Giroux 2007). As an educational approach "concerned with the intersection of power, identity, and knowledge" (Giroux 2007, p. 29), it advocates that all educators and teachers in a multicultural setting must be able to analyze the politics of education (Kincheloe and Steinberg 1997). As Freire puts it succinctly, critical pedagogy is about developing a critical consciousness by „learning to perceive social, political, and economic contradictions, and to take action against the oppressive elements of reality" (Freire 2000, p. 35). Given that societal power asymmetries can be (unknowingly) reproduced in educational settings, critical pedagogy argues for a humanizing classroom experience by advocating an educational philosophy that centers on the students' life experiences. Students are understood as "knowers" and "active participants," and the classroom should be an environment open for creativity, reflection, and action (Bartolome 1996). Through the assigned tasks, discussions, and self-reflection, educators and students are invited to develop a "critical ontology" and an "agency for change" of social and cultural issues (Meyer 2011; Kincheloe 2011).

Thus, the historical roots of Mathematics for Social Justice trace back to K-12 school education at the turn of the millennium. Later on, it became a topic of college classroom education. From the start, the goal was not to sacrifice mathematical content for social justice but to enrich the learning experience through critical pedagogy and to provide students with a positive mathematics identity (Buell and Shulman 2019). In college mathematics courses, projects have successfully been used to promote civic engagement among mathematics majors (Unfried and Canner 2019), while quantitative literacy and project-based general education courses have proven useful to equip students with the knowledge necessary to analyze fake news, misinformation, and promote good decision making in civil society (e.g., Branson 2019; O'Donovan and Geary 2019; Simic-Muller 2019). In calculus classes, students have shown particular responsiveness to case studies addressing social justice issues rather than to the classical physical problems that dominate much of calculus (Hoke et al. 2019). Common among all these approaches is the understanding that the student should be an active member of society, make moral decisions, use mathematics for social good, and be able to critically examine other peoples' mathematics while being able to reflect their own social background. For the non-STEM major, the goal is to have quantitative literate citizens. Robert Moses, one of the founding fathers of the Mathematics for Social Justice movement, understood mathematical literacy as a civil right (Moses and Cobb 2002). In turn, mathematics majors should become mathematicians with a social awareness of their work's cultural and social problems, including bias, discrimination, diversity, inclusion, fairness, equality, and equity.

Critical pedagogy itself traces its roots back to the critical theory developed by the Frankfurt school. Just like the writers from the Institute of Social Research at the University of Frankfurt perceived the world in an "urgent […] need of reinterpretation" after World War 1 and the economic collapse during the Weimar republic (Kincheloe and Mclaren 2011, p. 286), the recent writings from Mathematics for Social Justice perceive society (and thus the education system) in a similar need of reinterpretation and restructuring. In doing so, critical pedagogy does not build on the assumption of human objectivity and universal foundations for culture but generally rejects any such predisposition in its approach (Leistyna et al. 1996). Hence, it is unsurprising that Mathematics for Social Justice can also challenge existing beliefs about what counts as mathematics, including but not limited to its universality across cultures (Nolan and Graham 2021) and who can participate in its practice (Walshaw 2021). More generally, Mathematics for Social Justice is closely intertwined with questions of equity. So much so that the definitions of equity and social justice and the lines between them are often blurred or can even be non-existing for some authors (for an overview, we refer to Xenofontos et al. 2021).

To summarize, Mathematics for Social Justice is as much about "learning to play the game" as it is about "learning to change the game"(Gutiérrez 2009). It is a view towards education which is perhaps best described in the following fictional student's course evaluation for a social justice-themed mathematics course:

> "I finally get that […] teaching mathematics about, or through, social justice isn't just about poverty statistics and world population figures… it's also in the thoughts and actions of the teacher toward his/her students and in the thoughts and actions of students toward each other. It's about feeling safe to be who I am and, at the same time, to critically question who I want to become and what (and who) I value. And, most of all, I think it's also about opening up the content of mathematics […] to this same kind of critical questioning." (Nolan 2009, pp. 214–215)

### Paul Ernest's Philosophy of Mathematics Ethics Education

A similar view regarding the ethical challenges coming from mathematics can be found in the recent works of the philosopher Paul Ernest. In light of his perceived overvaluation of mathematics in society, Ernest calls for mathematicians to acknowledge the dark side of their subject and challenges the belief that practicing and learning more mathematics always leads to good outcomes (Ernest 2020b). Ernest (2016b) argues that the overt values ascribed to mathematics, including but not limited to truth, purity, universality, and beauty, can make it difficult for mathematicians to spot the covert values of objectism and ethics. In identifying three aspects of collateral damage to learning mathematics, he effectively bridges Ethics in Mathematics and Mathematics for Social Justice:

> "First, the nature of pure mathematics itself leads to styles of thinking that can be damaging when applied beyond mathematics to social and human issues. Second the applications of mathematics in society can be deleterious to our humanity unless very carefully monitored and checked. Third, the personal impact of learning mathematics on learners' thinking and life chances can be negative for a minority of less successful students, as well as potentially harmful for successful students." (Ernest 2018, p. 1)

For mathematics to have the appropriate valuation in society, we both need people trying to actively do good with mathematics and those trying to prevent harm. Therefore, it needs quantitatively literate, empowered, and self-reflected citizens who can question the existing uses of mathematics. But it also needs mathematicians who understand the power and dangers of their mathematics and applications and who are willing to consider the limitations of mathematical problem-solving approaches in social and human issues. This situates Paul Ernest in a close relationship to both Ethics in Mathematics and Mathematics for Social Justice.

## Comparison

Despite the commonly stated goal of fostering moral decision-making in students and mathematical practitioners, the assumptions that go into Ethics in Mathematics and Mathematics for Social Justice, their educational movements, and their teaching philosophies are quite different. In what follows, we outline some commonalities and differences and observe that their approaches to teaching and target groups can be understood as complementary instead of conflicting.

Ethics in Mathematics and Mathematics for Social Justice both attempt to establish an agency for change in their students. However, the message coming out of Mathematics for Social Justice is more positive than the one uttered by Ethics in Mathematics, focusing on enabling students to do good in society rather than teaching them how to prevent harm. It understands students as members of society who need to be empowered. In contrast, Ethics in Mathematics focuses on teaching those already in power: people with robust mathematical training are considered powerful members of our civic society who need to be taught to wield the "double-edged sword of mathematics" carefully (Chiodo 2/13/2020). The rather technocratic view of Ethics in Mathematics as ethics from mathematicians for mathematicians is complemented by the non-technocratic and empowering messages from Mathematics for Social Justice. Through both Ethics in Mathematics and Mathematics for Social Justice, students are trained to have a social awareness of their work's cultural and social problems, but the common-sense assumptions behind Ethics in Mathematics are philosophically much lighter than those going into the traditional Mathematics for Social Justice curricula. By participating in the Ethics in Mathematics seminars, students are supposed to level up in their ethical awareness and actions. However, the Ethics in Mathematics lectures only marginally address many of the workplace issues in the focus of Mathematics for Social Justice, including diversity, inclusivity, and equity. Hence, the Ethics in Mathematics lectures should at least be complemented by professional ethics courses dealing with these issues.

Mathematics for Social Justice is the larger scientific and educational movement, drawing from critical theory and pedagogy and the activist effort of people in the fight for social justice in politics and society. The growing literature on Ethics in Mathematics, on the contrary, is primarily being developed by practicing mathematical researchers. It does not challenge claims for universal foundations of reason and rational objectivity, and their work is mostly independent of foundational assumptions about the logical nature of mathematics. Unlike Paul Ernest's philosophy of mathematics education and many aspects of Mathematics for Social Justice, Ethics in Mathematics is only in minor parts social constructivist, in that it accepts that definitions, lemma, propositions, and theorems are coming from a research agenda designed by humans informed by their environment and culture. However, this is where their social constructivism ends, as the mathematical theorems of pure mathematics are still seen as universal within a specific framework of logic and set theory. Thus, the Ethics in Mathematics perspective on ethics of mathematics rather closely follows Gold's criticism of social constructivist philosophies of mathematics (Gold 1999). It tries to distinguish mathematics from mathematical knowledge and builds much on its ethical reasoning on the usefulness of mathematics. The human aspects are restricted to research, teaching and doing mathematics, and applying mathematics to extra-mathematical problems. The Platonist beliefs commonly held by mathematicians (Sriraman 2004) are not actively challenged, but definitory questions about the nature of pure mathematics have been raised (Müller 2018). The Ethics in Mathematics lectures only implicitly deal with Paul Ernest's view that different philosophies of mathematics, views ascribed to mathematics, and mathematics participation reinforce or hinder one another (Ernest 2003).

## Integration into the Curriculum

Only a few practical ethics of mathematics courses have been outlined in the literature (Franklin 2005; Chiodo and Bursill-Hall 2019), but the necessity to properly interweave them with the existing mathematical curriculum is stressed in the most recent works (Chiodo and Vyas 2019; Chiodo and Bursill-Hall 2019). Any course on ethics must also take into account the other existing possibilities through which students can learn about these issues. Do they have to take a general ethics course focusing on proper scientific conduct? Are there accompanying philosophy and history of science classes on which the ethics lectures can build? Do the students have to take a minor subject or fulfill general education requirements? If any of these exists, those designing the curriculum must be very careful about the likelihood of knowledge transfer. Are students likely to transfer their knowledge to the ethics of mathematics? The Ethics in Mathematics founders argued that this could be particularly difficult for students if the mathematics faculty does not fully support the ethics course or the ethical elements found in other courses (Chiodo and Bursill-Hall 2019).

Unlike the standard introductory courses on linear algebra and analysis that most mathematics students see at the beginning of their studies, introductory ethics courses have to be adjusted to the existing curriculum and the already existing (or non-existing) informal discourse on the ethics of mathematics found within the faculty. This can be difficult, and it also holds for attempts to introduce elements from Mathematics for Social Justice into the curriculum. Only recently, Mathematics for Social Justice has incited disagreement among parts of the mathematical community, with one commenting author calling it another "schism in the field" (Crowell 2022). These recent debates have only underlined the difference between introducing an individual course and interweaving its ideas with the whole curriculum. It is one of the core challenges that any ethics for mathematicians course faces. Ethics in Mathematics lectures at Cambridge found a way through a series of compromises: they are non-examinable, complemented by history of science and mathematics classes, and can build on a system of student societies found at few other universities. Being a seminar course specifically designed for and led by mathematicians, Ethics in Mathematics was able to reduce the alienation felt among the participating students.

## Conclusion

It is not surprising that ethical analyses of mathematics have regained prominence outside the educational sciences in recent years. The moral questions that mathematics has imposed on our society have led to new approaches, the most important of which are perhaps simple and old insights. Yet, it is not enough to look at individual subfields, let alone individual formulae, of mathematics when the digitalization and automatization of the economy and our everyday lives are driving its mathematization ever further.

The commonly held beliefs about good mathematics deeply influence what mathematicians understand as morally good mathematics and affect their agency for change which they feel towards ethical aspects of their field. Despite not formally taking on a specific philosophy of education, the Cambridge University Ethics in Mathematics Project, as ethics from mathematicians for mathematicians, nonetheless take on subtle positions which can be found in existing philosophies of mathematics educations and philosophies of mathematics. Coming out of practical considerations, these positions are closely related to the general mathematical practitioners' positions. Ethics in Mathematics and Mathematics for Social Justice are different perspectives on ethical responsibility in mathematics. Together they provide fruitful ground for discussions and shed light on the different social, cultural, and ethical issues surrounding mathematical practice and education. Our analysis outlined some fundamental differences and provided suggestions for the missing pieces in the Ethics in Mathematics research agenda. Ethics in Mathematics and Mathematics for Social Justice should be seen as complementary, using educational philosophies specifically selected to fit their target groups.

# Publication bibliography